\documentclass[12pt]{article}
\usepackage[dvips]{color} 
\usepackage{amsmath, amsthm, amssymb}
\usepackage[left=1.5in,top=1in,right=1in,nohead]{geometry}
\usepackage{hyperref}
\usepackage{latexsym}
\usepackage{mathrsfs}
\usepackage{rawfonts}
\usepackage[dvips,dvipdfmx]{graphicx}
\usepackage{float}
\usepackage[all]{xy}
\usepackage[sc,osf]{mathpazo}  

\usepackage{amsmath,amsthm,amssymb}
\usepackage{hyperref}
\usepackage{latexsym}
\usepackage{arydshln}

\usepackage{bmpsize}

\def\bct{\begin{center}}
\def\ect{\end{center}}
\def\bit{ \begin{itemize} }
\def\eit{ \end{itemize} }

\def\<{\langle}
\def\>{\rangle}

\newtheorem{theorem}{Theorem}
\theoremstyle{plain}

\newtheorem{definition}{Definition}

\numberwithin{equation}{section}

\title
{On K{\"a}hler Structures of Taub-NUT and Kerr Spaces}
\author{\" Ozg\"{u}r Kelek\c{c}i}
\date{ \textit{ {\small 
    Department of Basic Sciences \& Faculty of Engineering, \\ University of Turkish Aeronautical Association,  Ankara,  Turkey} \\
     {\small okelekci@thk.edu.tr}}
     }

\begin{document}
\maketitle

\begin{abstract}
In this paper, we study the K{\"a}hlerian nature of Taub-NUT and Kerr spaces which are gravitational instanton and black hole solutions in general relativity. We show that Euclidean Taub-NUT metric is hyper-K{\"a}hler with respect to the usual almost complex structures by employing an alternative explicit coframe, and Euclidean Kerr metric is globally conformally K{\"a}hler. We also show that conformally scaled Euclidean Kerr space admits a K{\"a}hler structure by applying a conformal scaling factor stemming from the Lee-form of the original metric or alternatively a factor coming from self-dual part of the Weyl tensor $(W^+)$.
\end{abstract}

 \noindent {\it Keywords}: Kerr space; Taub-NUT space; K{\"a}hler manifold; Hyper-K{\"a}hler manifold; locally conformally K{\"a}hler; Lee form.

\section{Introduction}

Taub-NUT and Kerr spaces, in addition to their significance in general relativity, possess rich geometric properties which provide interesting mathematical research problems. The Taub-NUT metrics were discovered by Taub \cite{a1} and later extended by Newmann-Unti-Tamburino \cite{a2}. They are complete Ricci-flat self-dual metrics on $\mathbb{R}^4$ and Hawking showed that they are non-trivial examples of the gravitational instantons \cite{a3}. Taub-NUT space has attracted a significant interest from both mathematics and physics research community \cite{a4,a5,a6,a7,a8}. Despite the existence of several studies on them there are few studies which are concentrated on their explicit hyper-K{\"a}hler structure. Gaeta etal. obtained Taub-NUT space from Euclidean eight-dimensional space $\mathbb{R}^8$  through a momentum map construction by implementing the HKLR theorem (Hitchin, Karlhede, Lindström, Roček) \cite{a7,a9}. However, there is no explicit connection of that study to the widely used Taub-NUT metrics such as given in \cite{a4} and \cite{a10}.

The Kerr metric mathematically describes a type of rotating black hole, and it generalizes the rotating version of Schwarzschild metric \cite{a11,a12}. Recently, geometric modeling of observed M87* black hole as a Kerr black hole was investigated \cite{KerrBH}.   Our focus will be on Euclidean versions of Kerr and Taub-NUT spaces rather than Lorentzian ones since we are searching for integrable almost complex structures. Lorentzian metrics can be transformed into Euclidean metrics by a coordinate transformation of the form  $t \rightarrow i \tau$, known as Wick rotation in physics literature \cite{a13}. Gibbons etal. discovered Kerr-Taub-bolt instanton solution as a rotating generalization of the Taub-bolt solution in a similar way that the Kerr solution generalizes the Schwarzschild solution \cite{a14}. Dixon showed that Wick-rotated Kerr metric, also called Riemannian Kerr metric, is ambitoric \cite{a15}. Aliev etal studied the general Kerr-Taub-bolt instanton and obtained a harmonic self-dual 2-form \cite{a16}. However, it is not clear whether the obtained self-dual 2-form is the genuine K{\"a}hler form with respect to the given metric and there is no mention of its almost complex structure. Babak etal studied some K{\"a}hler metrics via Lorentzian geometry where they have given K{\"a}hler structure only on some specific sections of Kerr space \cite{a17}. Despite the existence of these studies about K{\"a}hler structures of Kerr and Taub-NUT spaces, alternative explicit derivations of their K{\"a}hler structures starting directly from the commonly used metrics will be helpful for future studies.

In this paper, we first introduce an alternative explicit co-frame for the commonly used form of Euclidean Taub-NUT metric. We explicitly compute self-dual K{\"a}hler forms and corresponding almost complex structures $(J_1, J_2, J_3)$ by using this alternative co-frame. We also show that they are integrable and satisfy the quaternionic relations which completes the hyper-K{\"a}hler structure. Second, we focus on Wick-rotated Kerr metric (Euclidean Kerr) and show that it is globally conformally K{\"a}hler by calculating its K{\"a}hler form, Lee-form, and almost complex structure. Then, by applying a conformal scaling factor we show that there exists a genuine closed K{\"a}hler form and integrable almost complex structure for the conformally scaled Euclidean Kerr metric. As another verification we computed matrix form of $W^+$ and showed that the conformal factor obtained from exact Lee-form matches with the conformal factor of Derdzi\'{n}ski's theorem.

\section{Preliminaries}

In this section we shall give some definitions and theorems needed to prove our results. Note that we use Einstein summation convention throughout this paper, i.e. summing over repeated indices is assumed.

\begin{definition} Let M be an n-dimensional complex manifold and let ${z^\mu}$ be local coordinates on a coordinate neighborhood U. A mixed tensor ${J_\alpha}^\beta$ can be defined such that
\begin{equation}
    J=i\ d z^\mu \otimes\frac{\partial}{\partial z^\mu}-\ i\ d {\bar{z}}^\mu \otimes\frac{\partial}{\partial{\bar{z}}^\mu}  ~\ , \quad  \mu=1,2,...,n
\end{equation}
J is called almost complex structure and has the following properties: \begin{itemize}
\item[(i)] J is a tensor. 
\item[(ii)] J is real.
\item[(iii)] ${J_\alpha}^\beta\ {J_\beta}^\sigma=-{\delta_\alpha}^\sigma  ~\ , \quad     \alpha\ ,\ \sigma= 1,...\ ,2n.$
\end{itemize} 
\end{definition}

\begin{theorem} \label{thm1}
An almost complex structure J on a manifold $M$ is integrable if and only if the Nijenhuis tensor $N(X,Y)$ vanishes for any two vector fields X and Y \cite{a18}. Nijenhuis tensor is defined for any vector fields $X, Y\ \in T_pM$ as follows
 \begin{equation} \label{Nijenhuis}
        N\left(X,Y\right) = \left[X,Y\right]+J\left[JX,Y\right]+J\left[X,JY\right]-\left[JX,JY\right]  
     \end{equation}
\end{theorem}

\begin{theorem} \label{thm2}
 Let $(M, J)$ be an almost complex manifold. If J is integrable, the manifold $M$ is complex. This is also known as Newlander-Nirenberg theorem \cite{a19}.
\end{theorem}

\begin{definition} \label{defKahler}
Let  $\left(M,J\right)$ be a complex manifold. A Riemannian metric on M is called Hermitian if it is compatible with the complex structure J of M,

\begin{equation} 
g(JX,JY)=g(X,Y) \ ~ \text{for all} \ ~ X,Y \in T_p M \ ~ \text{and } \ ~  p\in M. 
\end{equation}
\end{definition}
\noindent Then the associated differential two-form $\omega $ defined by
\begin{equation} \label{omega2}
\omega(X,Y) =  g(JX,Y)
\end{equation} 
is called the {\em K\"{a}hler form} of the Hermitian metric. 
If moreover the form $\omega$ is closed, the metric is called {\em K\"ahler}.  

\begin{definition} \label{hyperKahler}
Let $(M,g)$ be  a 4k-dimensional Hermitian manifold admitting three almost complex structures I, J, K. $(M,g)$ is called  hyper-K\"ahler if I, J, K satisfies the quaternionic relations such as $I^2=J^2=K^2 =-Id, I.J = -J.I = K$. Moreover it can be shown that $(M,g)$ is hyper-K\"ahler if and only if $Hol(g) \subset Sp(k)$ \cite{a18}.

\end{definition}

\begin{theorem} \label{Lorentz-Hermitian}
A manifold with a metric of Lorentz signature cannot 
admit an almost Hermitian structure \cite{a20}. 
\end{theorem}

\begin{theorem} \label{Lckthm}
 The Hermitian manifold $(M^{2n}, J, g)$  is locally conformally K{\"a}hler (l.c.K) if and only if there exists a globally defined closed 1-form (called Lee form)  $\xi$ on $M^{2n}$ so that  \cite{a21}
\end{theorem}
\begin{equation}\label{Lee-form1}
 d\ \omega=\xi\wedge\omega  
\end{equation}
                                                
\noindent Furthermore, the Lee form $\xi$ is uniquely determined by any of the following formulae \cite{a22}.
\begin{equation} \label{Lee-form2}                  \xi=-\frac{1}{n-1}\left(\delta\omega\right)\circ J \quad \text{or} \quad \xi_i=-\frac{2}{m-2}{(\nabla}_{\alpha\ }J_\beta^\alpha){\ J}_i^\beta  
\end{equation}

\noindent where $\delta$ is the co-differential operator and $\text{dim}\left(M\right)=m=2n$.  If moreover the Lee form $\xi$ is exact then $(M^{2n}, J, g)$ is globally conformally K{\"a}hler. On a l.c.K. manifold, differentials of the negative of the logarithm of local conformal factors can be shown to agree on the overlaps so to form a global 1-form $\xi$. Conversely, starting with a closed 1-form on $M^{2n}$ satisfying \eqref{Lee-form1}, Poincaré lemma guarantees that there is an open cover $\mathcal{U}=\{U_k\}_{k\in I}$ of  $M^{2n}$ and a family of $C^\infty$ functions $ f_k:U_k\rightarrow\mathbb{R}$ so that $\xi=df_k$ on $U_k$. Also, by multiplying $\omega$ with $e^{-f_k}$ one gets $d(e^{-f_k}\ \omega)=0$ and hence $\hat{g}=e^{-f_k}\ g$ becomes a K{\"a}hler metric on $U_k$. 

\begin{theorem} \label{Derdzinski}
Let $(M,g)$ be an oriented 4-dimensional 
Einstein manifold. Assume that the self-dual Weyl tensor $W^+$ does not vanish. Then 
there exist two equivalent properties   \cite{Derz}

i) $W^+$ has, at each point, at most two distinct eigenvalues on $\Lambda_2 ^+$ 

ii) up to a double cover, the metric $|W^+|^{2/3}g$  is K{\"a}hler (with respect to some complex structure on $M$ compatible with the original orientation). 

\end{theorem}

where $\Lambda_2 ^+$ is bundle of self-dual 2-forms on $M$.

\noindent Theorem~\ref{Derdzinski} can be used for an alternative verification of the conformal K{\"a}hler structure. Let us define self-dual ($\Omega_+$) and anti-self-dual ($\Omega_-$) bases   for $\Lambda^2$ (bundle of 2-forms on $M$) using orthonormal basis one-forms $e^i$ as,

\begin{equation*}\label{w2basis}
\Omega_\pm=\left\{e^1 \wedge e^2 \pm e^3 \wedge e^4, e^1 \wedge e^3 \pm e^4 \wedge e^2, e^1 \wedge e^4 \pm e^2 \wedge e^3  \right\}
\end{equation*}

\noindent The Riemann curvature tensor defines a self-adjoint transformation $\mathcal{R}: \Lambda^2 \rightarrow \Lambda^2$ and it can be considered as a $6 \times 6$ block matrix with respect to the decomposition $\Lambda^2 = \Lambda_+^2 \oplus \Lambda_-^2$ \cite{AHS}

\begin{equation}
  \mathcal{R} \Lambda^2  =\biggl(  \begin{array}{c;{2pt/2pt}c} 
  W^+ + \frac{R}{12} I \, \, & 
 \mathring{Ric} \\ 
  \hdashline[2pt/2pt]
   \overset{}{\mathring{Ric}} &  W^- + \frac{R}{12} I 
\end{array} \biggr) \biggl( \begin{array}{c} \Omega_+ \\ 
  \hdashline[2pt/2pt] \Omega_- \end{array} \biggr)
\end{equation}

\noindent where $W^\pm$ are the self-dual and anti-self-dual parts of Weyl tensor in matrix form, $I$ is the identity map of $\Lambda^2$ and $\mathring{Ric}$ is the trace-free Ricci curvature ($Ric-\frac{R}{4}g$). Upper left block matrix is also represented by a $3 \times 3$ matrix $A$ whose entries can be computed by the following equation \cite{a10}
\begin{eqnarray}\label{Weyl+}
    A_{ij}&=&(R_{0i0j}+\frac{1}{2}\epsilon_{jkl} R_{0ikl})+\frac{1}{2}\epsilon_{imn}(R_{mn0j}+ \epsilon_{jkl} R_{mnkl})  \\
     A&=& W^+ + \frac{R}{12} I \nonumber
\end{eqnarray}

\section{Main Results}

\subsection{Hyper-K{\"a}hler Structure Of Euclidean Taub-Nut Metric} 

We will study Euclidean Taub-NUT space and use the metric given in \cite{a4}. Euclidean Taub-NUT metric is given in coordinates as 

\begin{equation}\label{TNut1}
\tilde{g}=V\left(dx^2+dy^2+dz^2\right)+V^{-1}\left(dt+\Theta\right)^2
\end{equation}
\begin{equation}\label{TNut2}
g_{TN}=\frac{\rho+1}{4\rho}d\rho^2+\rho\left(\rho+1\right)\left(\sigma_1^2+\sigma_2^2\right)+\frac{\rho}{\rho+1}\sigma_3^2                                    
\end{equation}

\noindent where  $V=1+1/2r$  and $\Theta$ is a 1-form on $\mathbb{R}^3$   satisfying $d \Theta=\ast \ dV$. x, y  and z are coordinates parameterizing $ \mathbb{R}^3 $  while t being another vertical coordinate. $ r $  is the radial coordinate for  $\mathbb{R}^3$ satisfying $r^2=x^2+y^2+z^2$. It is easy to show that the metrics in \eqref{TNut1} and \eqref{TNut2} are isometric by using the transformation map $f:\left(x,y,z,t\right)\mapsto\left(\rho,\theta,\phi,\psi\right)$ defined by the following relations
\begin{eqnarray}
x&=&\frac{\rho}{2}\sin{\theta}\cos{\phi} \ , \ \ y=\frac{\rho}{2}\sin{\theta}\sin{\phi} \ ,\ \ z=\frac{\rho}{2}\cos{\theta}, \nonumber \\ 
V&=&1+\frac{1}{\rho}\ ,\ \ \Theta=\frac{1}{2}\cos{\theta}d\phi\ , \ \ t=\psi/2 
\end{eqnarray}

\noindent $\sigma_1,\sigma_2$ and $\sigma_3$ given in \eqref{TNut2} are left invariant one-forms parameterizing the spaces with SU(2)  isometry and they are written in Euler angles $(\theta,\phi,\psi)$ as
\begin{eqnarray}
\sigma _{1} &=&(\sin \psi \,d\theta -\sin \theta \cos \psi \,d\phi )/2 \\ \nonumber
\sigma _{2} &=&(\cos \psi \,d\theta +\sin \theta \sin \psi \,d\phi )/2 \\ 
\sigma _{3} &=&(d\psi +\cos \theta \,d\phi )/2\nonumber
\end{eqnarray}

\noindent These one-forms satisfy the structure equations $d\sigma_i=\epsilon_{ijk} \ \ \sigma_j\wedge\sigma_k$ , where $\epsilon_{ijk}$ is the 3-dimensional anti-symmetric tensor. Euclidean Taub-NUT metric is also given in \cite{a10} following \cite{a3} as 

\begin{equation} \label{TNut3}
 \frac{r+m}{4\left(r-m\right)}{\rm dr}^2+\left(r^2-m^2\right)\left({\sigma_1}^2+{\sigma_2}^2\right)+4m^2\frac{r-m}{r+m}{\sigma_3}^2       
\end{equation}

\noindent The metric in \eqref{TNut2} can be obtained by setting $r=\rho+m$ and $m=1/2$ in \eqref{TNut3}. Now that we have obtained isometries between the most common explicit Taub-NUT metrics in literature we will continue with \eqref{TNut2}. As pointed out in \cite{a4} the 2-form $\omega$ obtained with respect to the most obvious integrable almost complex structure $\left(J\sigma_1\rightarrow\sigma_2,J\sigma_3\rightarrow-d\rho\left(1+\rho\right)/2 \rho \right)$ is not closed. Therefore, one must find a suitable co-frame for the metric and compatible almost complex structure so that K{\"a}hler conditions are satisfied. We write the following orthonormal co-frame in explicit coordinates as follows

\begin{eqnarray} \label{TNframe}
 e^1&=&\sqrt{\frac{\rho+1}{4\rho}}\left(\rho\cos{\theta}\cos{\phi}d\theta+\sin{\theta}\left(\cos{\phi}d\rho-\rho\sin{\phi}d\phi\right)\right)  \\ \nonumber
e^2&=&\sqrt{\frac{\rho+1}{4\rho}}\left(\rho\cos{\theta}\sin{\phi}d\theta+\sin{\theta}\left(\sin{\phi}d\rho+\rho\cos{\phi}d\phi\right)\right)  \\ \nonumber
e^3&=&\sqrt{\frac{\rho+1}{4\rho}}\left(-\rho\sin{\theta}d\theta+\cos{\theta}d\rho\right)   \ , \quad e^4=\sqrt{\frac{\rho}{4\left(\rho+1\right)}}\left(\cos{\theta}d\phi+d\psi\right)       
\end{eqnarray}

\noindent The vector fields dual to the co-frame $\{e^1,e^2,e^3,e^4\}$ are obtained as 

\begin{eqnarray}
e_1&=&\frac{2}{\sqrt{\rho\left(1+\rho\right)}}\left(\cos\theta \ \cos\phi{\ \partial}_\theta + \rho \ \sin\theta \ \cos\phi \ \partial_\rho+ \sin\phi \  \cot\theta \ \partial_\psi-\frac{\sin\phi}{\sin\theta}\partial_\phi\right) \nonumber \\ 
e_2&=&\frac{2}{\sqrt{\rho\left(1+\rho\right)}}\left(\cos\theta \ \sin\phi{\ \partial}_\theta+\rho \ \sin\theta \ \sin\phi \ \partial_\rho-\cos\phi \  \cot\theta \  \partial_\psi+\frac{\cos\phi}{\sin\theta}\partial_\phi\right) \nonumber \\ 
e_3&=&\frac{2}{\sqrt{\rho\left(1+\rho \right)}}\left(\rho \ \cos\theta{\ \partial}_\rho-\sin\theta \ \partial_\theta \right)  \ , \quad e_4=2\sqrt{\frac{\rho+1}{\rho}} \ \partial_\psi    
\end{eqnarray}

\noindent We choose the orientation $(e^1\wedge e^2\wedge e^3\wedge e^4)$ so that the following 2-forms become self-dual.
\begin{equation}\label{TNomega}
    \omega_1=e^1\wedge e^2+e^3\wedge e^4 \ , \ {\ \omega}_2=e^1\wedge e^4+e^2\wedge e^3 \ , \ {\ \omega}_3=e^1\wedge e^3+e^4\wedge e^2
\end{equation}

\noindent Almost complex structures that will yield the 2-forms given in \eqref{TNomega} with respect to the chosen orientation will be
\begin{eqnarray}\label{Jframe}
&{\hat{J}}_1 e_1=e_2 \ , \ {\hat{J}}_1 e_2=-e_1 \ , \    {\hat{J}}_1e_3=e_4  \ , \  {\hat{J}}_1e_4=-e_3 \\ \nonumber
&{\hat{J}}_2 e_1=e_4   \ , \  {\hat{J}}_2 e_2=e_3  \ , \ \ \ {\hat{J}}_2 e_3={-e}_2  \ , \  {\hat{J}}_2 e_4={-e}_1 \\ \nonumber
&{\hat{J}}_3 e_1=e_3   \ , \ {\hat{J}}_3 e_2={-e}_4   \ , \ {\hat{J}}_3 e_3={-e}_1  \ , \ {\hat{J}}_3 e_4=e_2 
\end{eqnarray}

\noindent Self-dual 2-forms given in \eqref{TNomega} can be written in coordinates
\noindent \begin{equation}\label{TNomega2}
\begin{split}
&\omega_1=\frac{1}{4}\biggl(\cos\theta \ d\rho \wedge d\psi+\rho  \sin\theta\left(\rho \cos\theta \ d\theta \wedge d\phi-d\theta \wedge  d\psi+ \sin\theta \ d\rho \wedge d\phi\right)+d\rho \wedge  d\phi \biggr)  \\ 
&\omega_2=\frac{1}{4}\biggl[\left(1+\rho\right) \sin\phi \ d\theta \wedge d\rho+\rho \cos\phi (1+\rho\ \sin^2\theta)\ d\theta \wedge d\phi+\rho \cos\phi  \cos\theta\left(d\theta \wedge d\psi- \sin\theta \ d\rho \wedge d\phi\right)  \\ 
& \quad \ \ \ \ + \cos\phi  \sin\theta \ d\rho \wedge d\psi-\rho \sin\phi  \sin\theta \ d\phi \wedge d\psi \biggr]  \\ 
&\omega_3=\frac{1}{4}\biggl[(1+\rho) \cos\phi\ d\theta \wedge d\rho-\rho\  \sin\phi(1+\rho\ \sin^2\theta)\ d\theta \wedge d\phi+\rho  \sin\phi\  \cos\theta \ ( \sin\theta\ d\rho \wedge d\phi-d\theta \wedge d\psi)\  \\  
& \quad \ \ -  \sin\phi\  \sin\theta\ d\rho \wedge d\psi)-\rho\  \cos\phi\  \sin\theta\ d\phi \wedge d\psi \biggr]   \end{split}            
\end{equation}                         

\noindent It is a straightforward computation to check that $\omega_1$, $\omega_2$ and $\omega_3$ are closed. Almost complex structures are given in \eqref{Jframe} with respect to co-frame basis elements. They should be written with respect to explicit coordinates as well. Following index equation can be written from definition of K{\"a}hler form.
\begin{eqnarray}\label{wcoeff}
\omega=\frac{1}{2}g_{\mu\nu}\ J_\sigma^\mu\ {\rm dx}^\sigma\wedge {\rm dx}^\nu=\frac{1}{2}\omega_{\sigma\nu}\ {\rm dx}^\sigma\wedge {\rm dx}^\nu \\ \nonumber
g^{\nu\alpha}g_{\mu\nu}\  J_\sigma^\mu=g^{\nu\alpha}\omega_{\sigma\nu}\ \Rightarrow\ \delta_\mu^\alpha{\ J}_\sigma^\mu={\ J}_\sigma^\alpha=g^{\nu\alpha}\omega_{\sigma\nu}
\end{eqnarray}

\noindent Coefficients of the $\omega_i$’s  ($\omega_{\sigma\nu}$) can be easily read from \eqref{TNomega2}. Then inverse metric coefficients ($g^{\nu\alpha}$) can be applied  in \eqref{wcoeff} to obtain the almost complex structures in coordinate frame. Coefficients of the J tensors can be expressed in matrix form,

\begin{eqnarray} \label{Jcoord1}
J_1&=&\left(\begin{matrix}0&\frac{1}{\rho}&0& \cos\theta\\-\frac{\rho\left(1+\rho\ \sin^2\theta\right)}{1+\rho}&0&-\frac{\rho\  \cos\theta  \sin\theta}{1+\rho}&0\\0& \cot\theta&0&-\frac{1+\rho{\sin}^2\theta}{ \sin\theta}\\-\frac{\rho \cos\theta}{1+\rho}&0&\frac{ \sin\theta}{1+\rho}&0\\\end{matrix}\right) \\ \nonumber
J_2&=&\left(\begin{matrix}0&-\frac{ \cos\phi \ \cot\theta}{\rho}&-\frac{ \sin\phi}{\rho}&\frac{(1+\rho\ \sin^2\theta) \cos\phi}{\rho \ \sin\theta}\\\frac{\rho^2 \sin2\theta  \cos\phi}{2(1+\rho)}&\frac{ \sin\phi \ \cot\theta}{1+\rho}&-\frac{(1+\rho\ \sin^2\theta) \cos\phi}{1+\rho}&-\frac{\left(\rho (2+\rho) \sin^2\theta + 1 \right) \sin\phi}{(1+\rho) \sin\theta}\\\rho \ \sin\phi& \cos\phi&0&\rho \ \cos\theta  \cos\phi\\-\frac{\rho \ \cos\phi  \sin\theta}{1+\rho}&\frac{ \sin\phi}{(1+\rho) \sin\theta}&-\frac{ \cos\theta  \cos\phi}{1+\rho}&-\frac{\cot\theta  \sin\phi}{1+\rho}\\\end{matrix} \right)  \\ \nonumber
J_3&=&\left(\begin{matrix}0&\frac{ \sin\phi \ \cot\theta}{\rho}&-\frac{ \cos\phi}{\rho}& -\frac{(1+\rho\ \sin^2\theta) \sin\phi}{\rho \ \sin\theta}\\-\frac{\rho^2 \sin2\theta  \sin\phi}{2(1+\rho)}&\frac{ \cos\phi \ \cot\theta}{1+\rho}& \frac{(1+\rho\ \sin^2\theta) \sin\phi}{1+\rho}&-\frac{\left(\rho (2+\rho) \sin^2\theta + 1 \right) \cos\phi}{(1+\rho) \sin\theta}\\\rho \ \cos\phi& -\sin\phi& 0 &-\rho \ \cos\theta  \sin\phi\\ \frac{\rho \ \sin\phi  \sin\theta}{1+\rho}&\frac{ \cos\phi}{(1+\rho) \sin\theta}& \frac{ \cos\theta  \sin\phi}{1+\rho}&-\frac{\cot\theta  \cos\phi}{1+\rho}\\\end{matrix} \right)
\end{eqnarray}

\noindent These almost complex structures satisfy the quaternionic relations such as ${J_1}^2={J_2}^2={J_3}^2=-\text{Id}$ , $J_1J_2={-J}_2J_1=J_3$ (others are cyclic). Lie brackets of basis vector fields $\{e_1,e_2,e_3,e_4\}$ are needed for Nijenhuis tensor calculations.
\noindent \begin{multline*}
\left[e_1,e_2\right]=\frac{2\left(\partial_\phi+\left(1+2\rho\right)\cos\theta \ \partial_\psi\right)}{\rho\left(1+\rho\right)^2}      \  , \ \ \left[e_1,e_4\right]=-\frac{2\ \cos\phi\ \sin\theta\ \partial_\psi}{\rho (1+\rho)}   
\hfill\end{multline*}~~\vspace{-6mm}
\begin{multline*}
\left[e_2,e_4\right]=-\frac{2\ \sin\phi\ \sin\theta\ \partial_\psi}{\rho (1+\rho)} \ , \  \left[e_3,e_4\right]=-\frac{2 \cos\theta \ \partial_\psi}{\rho (1+\rho)} \hfill\end{multline*}~~\vspace{-6mm}
\begin{multline*} 
\left[e_1,e_3\right]=\frac{2\cot\theta\ \sin\phi}{\rho{(1+\rho)}^2}\partial_\phi-\frac{2\cos\phi}{\rho{(1+\rho)}^2}\partial_\theta-\frac{2\sin \phi(\sin^2\theta(1+2\rho)+1)}{\rho \sin\theta{(1+\rho)}^2}\partial_\psi
\hfill\end{multline*}
\begin{multline} \label{brackets1}
\left[e_2,e_3\right]=-\frac{2\cot\theta\ \cos\phi}{\rho{(1+\rho)}^2}\partial_\phi-\frac{2\sin\phi}{\rho{(1+\rho)}^2}\partial_\theta+\frac{2\cos\phi(\sin^2\theta(1+2\rho)+1)}{\rho \sin\theta{(1+\rho)}^2}\partial_\psi
\hfill\end{multline}

\noindent Integrability condition for J requires Nijenhuis tensor to vanish for all vector fields as mentioned in Theorem~\ref{thm1}. It is now straightforward to show that Nijenhuis tensor is vanishing for all three almost complex structures by using the coordinate coefficients of $J_1$, $J_2$  and $J_3$ tensors given in \eqref{Jcoord1} and the brackets in \eqref{brackets1}. We show one example calculation for $J_1$, others can be computed similarly.

\noindent \begin{multline}
 N_1\left(e_1,e_3\right) = \left[e_1,e_3\right]+J_1 \left[J_1 \  e_1,e_3\right]+J_1 \left[e_1,J_1 \  e_3\right]-\left[J_1 \  e_1 ,J_1 \ e_3\right] \quad \ \ \
\hfill\end{multline}~~\vspace{-6mm}
\begin{multline*}
\phantom{xxxxxxxx} = \left[e_1,e_3\right]+J_1 \left[ e_2,e_3\right]+J_1 \left[e_1, e_4\right]-\left[e_2 ,e_4\right]
\hfill\end{multline*}~~\vspace{-6mm}
\begin{multline*}
\phantom{xxxxxxxx} =\left[e_1,e_3\right] -\frac{\cos\phi\ \sin2\theta}{{(1+\rho)}^2}\partial_\rho-\frac{2\cot\theta\ \sin\phi}{\rho{(1+\rho)}^2}\partial_\phi+\frac{2\cos\phi(1+\sin^2\theta)}{\rho{(1+\rho)}^2}\partial_\theta
\hfill\end{multline*}
\begin{multline*}
\phantom{xxxxxxxxxx}+\frac{2\sin\phi(1+\rho\ \sin^2\theta)}{\rho \sin\theta{(1+\rho)}^2}\partial_\psi + \frac{\cos\phi\ \sin2\theta}{{(1+\rho)}^2}\partial_\rho-\frac{2\cos\phi\ \sin^2\theta}{\rho{(1+\rho)}^2}\partial_\theta + \frac{2\sin\phi\ \sin^2\theta\ (1+\rho)}{\rho\ \sin\theta  {(1+\rho)}^2}\partial_\psi
\hfill\end{multline*}
\begin{multline*}
\phantom{xxxxxxxx}=\left[e_1,e_3\right]-\frac{2\cot\theta\ \sin\phi}{\rho{(1+\rho)}^2}\partial_\phi+\frac{2\cos\phi}{\rho{(1+\rho)}^2}\partial_\theta+\frac{2\sin \phi(\sin^2\theta(1+2\rho)+1)}{\rho \sin\theta{(1+\rho)}^2}\partial_\psi = 0.
\hfill\end{multline*}

\noindent Hence, this completes all the necessary conditions for hyper-K\"{a}hler structure, ensuring that the explicit Euclidean Taub-NUT metric \eqref{TNut2} is hyper-K\"{a}hler with respect to the complex structures obtained. Note that an alternative way to check the integrability is to show that J is covariantly constant $(\nabla J=0)$ with respect to the Levi-Civita connection associated with the metric.

\subsection{K{\"a}hler Structure Of Euclidean Kerr Metric} \, 

\noindent Kerr metric can be written using Boyer-Lindquist coordinates in the following form \cite{a23},
\begin{equation}\label{Kerr1}
\tilde{g}_K= \tilde{\Xi}\ \left(\frac{{\rm dr}^2}{\tilde{\Delta}}+{d\theta}^2\right)+\frac{\sin^2\theta}{\tilde{\Xi}}\left(-\alpha dt+\left(r^2+\alpha^2\right)d\phi\right)^2+\frac{\tilde{\Delta}}{\tilde{\Xi}}\left(dt-\alpha\sin^2\theta \ d\phi\right)^2
\end{equation}

\noindent where $\tilde{\Delta}=r^2-2Mr+\alpha^2$ and $\tilde{\Xi}=r^2+\alpha^2\cos^2{\theta}$. It describes the exterior gravitational field of a rotating black hole with mass M and angular momentum per unit mass of  $\alpha$ which are positive real parameters. Kerr space is a Lorentzian 4-manifold categorized under Petrov Type D class \cite{a23}. A manifold having a metric of Lorentz signature cannot admit an almost Hermitian structure as mentioned in Theorem~\ref{Lorentz-Hermitian}. In this case, one can apply Wick rotation to the metric and obtain another version of the given metric with Riemannian signature. As a result of this transformation, solutions of the Lorentzian Einstein equations are mapped into solutions of the Riemannian Einstein equations. Following \cite{a14} we apply Wick rotation to \eqref{Kerr1} and get,

\begin{equation} \label{Kerr2}
g_K=\Xi\left(\frac{{\rm dr}^2}{\Delta}+{d\theta}^2\right)+\frac{\sin^2\theta}{\Xi}\left(\alpha dt+\left(r^2-\alpha^2\right)d\phi\right)^2+\frac{\Delta}{\Xi}\left(dt-\alpha\sin^2\theta d\phi\right)^2
\end{equation}

\noindent Here Wick rotation is applied to both t and $\alpha$ such that $t \rightarrow it$ , $\alpha \rightarrow i\alpha$ to yield \eqref{Kerr2} with $\Delta=r^2-2Mr-\alpha^2$ and $\Xi=r^2-\alpha^2\cos^2{\theta}$. This metric is known as Euclidean (Riemannian) Kerr metric. Now we can search for Kähler forms and integrable almost complex structures. An orthonormal co-frame for this metric can be written as,
\begin{equation}
\{{e}^1,e^2,e^3,e^4\}=\left\{\sqrt{\frac{\Xi}{\Delta}}dr,\sqrt\Xi\ d\theta,\frac{\sin\theta}{\sqrt\Xi}\left(\alpha dt+\left(r^2{-\alpha}^2\right)d\phi\right),\sqrt{\frac{\Delta}{\Xi}}\left(dt-\alpha{\sin}^2\theta\ d\phi\right)\right\}         
\end{equation}

\noindent The vector fields dual to the co-frame basis elements are  
\begin{equation}
\{e_1,e_2,e_3,e_4\}=\left\{\sqrt{\frac{\Delta}{\Xi}} \ \partial_r,\frac{\partial_\theta}{\sqrt\Xi},\frac{\alpha\ {\sin}^2\theta \ \partial_t+\partial_\phi}{\sin\theta\ \sqrt\Xi},\frac{\left(-\alpha \ \partial_\phi+\left(r^2-\alpha^2\right)\partial_t\right)}{\sqrt{\Delta\ \Xi}}\right\}     
\end{equation}
  
\noindent We choose the orientation ($e^1\wedge e^2\wedge e^3\wedge e^4$) so that the 2-form $\omega=e^1\wedge e^4+e^2\wedge e^3$ become self-dual as we did in Taub-NUT space. And the corresponding almost complex structure is the following

\begin{equation} \label{Jframe2}
\hat{J}e_1=e_4 \ , \quad  \hat{J}e_2=e_3 \ , \quad  \hat{J}e_3={-e}_2 \ , \quad  \hat{J}e_4={-e}_1
\end{equation}

\noindent Self-dual 2-form $\omega$ and coefficients of almost complex structure  $\hat{J}$ can be written in coordinates,
\begin{equation}\label{wKerr1}
 \omega=dr\wedge dt-\alpha{\sin}^2\theta \ dr\wedge d\phi-\alpha \sin\theta \ dt\wedge d\theta+\left(r^2-\alpha^2\right)\sin\theta \ d\theta\wedge d\phi    
\end{equation}
\begin{equation}\label{JKerr}
J=\left(
\begin{array}{cccc}
 0 & -\frac{\Delta }{\Xi } & -\frac{\alpha  \sin \theta }{\Xi } & 0 \\
 \frac{r^2-\alpha ^2}{\Delta } & 0 & 0 & -\frac{\alpha }{\Delta } \\
 \alpha  \sin \theta & 0 & 0 &\csc \theta \\
 0 & \frac{\alpha  \Delta  \sin^2 \theta}{\Xi } & \frac{\sin \theta \left(\alpha ^2-r^2\right)}{\Xi } & 0 \\
\end{array}
\right)  
\end{equation}

\noindent Taking the exterior derivative of $\omega$ we get  $d\omega=2\left(r+\alpha\ \cos\theta\right)\sin\theta\ dr\wedge d\theta\wedge d\phi$, which is not identically vanishing so the metric under consideration is not K{\"a}hler. Indeed, there exists a closed self-dual 2-form defined on the Euclidean Kerr space as shown in \cite{a16}. We compute the corresponding almost complex structure for that closed 2-form as follows

\begin{equation}
\tilde{J}=\left(
\begin{array}{cccc}
0&-\frac{\Delta}{\Xi{(r-\alpha \cos\theta)}^2}&-\frac{\alpha \sin\theta}{\Xi{(r-\alpha \cos\theta)}^2}&0\\\frac{r^2-\alpha^2}{\Delta{(r-\alpha \cos\theta)}^2}&0&0&-\frac{\alpha}{\Delta{(r-\alpha \cos\theta)}^2}\\\frac{\alpha sin\theta}{{(r-\alpha \cos\theta)}^2}&0&0&\frac{csc\theta}{{(r-\alpha \cos\theta)}^2}\\0&\frac{\alpha\Delta \sin^2 \theta}{\Xi{(r-\alpha \cos\theta)}^2}&\frac{\sin \theta (\alpha^2-r^2)}{\Xi{(r-\alpha \cos\theta)}^2}&0\\
\end{array}
\right)
\end{equation}

\noindent It turns out that $\tilde{J}$ is only  a scaled version of \eqref{JKerr}. However, $\tilde{J}$  does not even satisfy the condition $\tilde{J}^2=-\text{Id}$ for an almost complex structure, let alone the integrability condition. Therefore, Euclidean Kerr metric fails to exhibit genuine K{\"a}hlerian nature in its current form given in \eqref{Kerr2}. 

\noindent Theorem~\ref{Lckthm} is useful for doing further analysis on K{\"a}hlerity of a given space. The Lee form associated with K{\"a}hler form $\omega$ in \eqref{wKerr1} is calculated by using the formula in \eqref{Lee-form2} as
\begin{equation} \label{LeeKerr}
\xi=\frac{2\left(dr+\alpha\ \sin\theta\ d\theta\right)}{r-\alpha\ \cos\theta}=d\ \left(\text{Log}\left[ C \left(r-\alpha \cos\theta\right)^2\right]\right)
\end{equation}

\noindent where $C$ is an arbitrary constant. The existence of an exact Lee form $\xi$ ensures that the Euclidean Kerr metric given in \eqref{Kerr2} is globally conformally K{\"a}hler due to Theorem~\ref{Lckthm}. Note that the obtained Lee form also satisfies \eqref{Lee-form1}. Taking $C=1$ \footnote{arbitrary constant $C$ could be kept in the proceeding calculation steps but it doesn't affect the obtained results.}, conformal scaling factor can be determined from \eqref{LeeKerr}  as

\begin{equation}\label{confactor}
\lambda=exp(-\text{Log}\left[\left(r-\alpha \cos\theta\right)^2\right])={1/\left(r-\alpha \cos\theta\right)}^2
\end{equation}

\noindent When this conformal factor is applied to the metric in \eqref{Kerr2}  the resulting metric is expected to be K{\"a}hler in accordance with the Theorem~\ref{Lckthm} and its sketch of proof given thereafter. Applying this conformal factor, we get the scaled metric $ g_{conf}= g_K/\left(r-\alpha \cos\theta \right)^2$, corresponding co-frame and vector frame as follows
\noindent \begin{equation*}
 \{e^1,e^2,e^3,e^4\}=\frac{1}{\left(r-\alpha \cos\theta\right)}\left\{\sqrt{\frac{\Xi}{\Delta}}dr,\sqrt\Xi\ d\theta,\frac{\sin\theta}{\sqrt\Xi}\left(\alpha dt+\left(r^2{-\alpha}^2\right)d\phi\right),\sqrt{\frac{\Delta}{\Xi}}\left(dt-\alpha{\sin}^2\theta\ d\phi\right)\right\} \end{equation*}

\begin{equation}
\{e_1,e_2,e_3,e_4\}=\left(r-\alpha \cos\theta\right)\left\{\sqrt{\frac{\Delta}{\Xi}} \ \partial_r,\frac{\partial_\theta}{\sqrt\Xi},\frac{\alpha\ {\sin}^2\theta \ \partial_t+\partial_\phi}{\sin\theta\ \sqrt\Xi},\frac{\left(-\alpha \ \partial_\phi+\left(r^2-\alpha^2\right)\partial_t\right)}{\sqrt{\Delta\ \Xi}}\right\}
\end{equation}

\noindent Same orientation and almost complex structure in \eqref{Jframe2} will be used to compute $\hat{\omega}$ in coordinate frame as

\begin{equation}\label{wKerr2}
\hat{\omega}= \frac{1}{\left(r-\alpha \cos\theta\right)^2} \left( dr\wedge dt-\alpha{\sin}^2\theta \ dr\wedge d\phi-\alpha \sin\theta \ dt\wedge d\theta+\left(r^2-\alpha^2\right)\sin\theta \ d\theta\wedge d\phi \right)   
\end{equation}

\noindent One can easily verify that $\hat{\omega}$ is now closed. This is already expected from the construction of new metric as explained in sketch of proof of Theorem~\ref{Lckthm}. Note that the coordinate-based coefficients of almost complex structure in \eqref{JKerr} does not change due to having the same factor for each element of the scaled vector frame. We compute the Lie brackets of the basis vector fields,

\noindent \begin{multline*}
\left[e_1,e_2\right]=\frac{\alpha \sqrt{\Delta}\ (\cos\theta\ \partial_\theta -r \sin\theta\ \partial_r)}{{(r+\alpha\ \cos\theta)}^2} \quad , \quad  \left[e_2,e_4\right]=\frac{r\ \alpha\ \sin\theta \left((r^2-\alpha^2)\ \partial_t-\alpha\ \partial_\phi \right)}{\sqrt{\Delta} \ {(r+\alpha\ \cos\theta)}^2}
\hfill\end{multline*}
\begin{multline*}
\left[e_1,e_3\right]=\frac{\alpha \sqrt{\Delta}\ \cot \theta\ (\partial_\phi + \alpha \sin^2 \theta\ \partial_t)}{{(r+\alpha\ \cos\theta)}^2} \quad , \quad  \left[e_3,e_4\right]=0
\hfill\end{multline*}
\begin{multline*}
\left[e_1,e_4\right]=-\frac{\alpha^2 \cos\theta\ \Delta\ +\alpha\ \Xi\ (M-r)}{\Delta{(r+\alpha\ \cos\theta)}^2}\partial_\phi+\frac{\alpha\ \cos\theta\ \Delta\ (r^2-\alpha^2)+\Xi(\Delta r-M(r^2+\alpha^2))}{\Delta{(r+\alpha  \cos\theta)}^2}\partial_t
\hfill\end{multline*} 
\begin{multline}
 \left[e_2,e_3\right]= \frac{\alpha\ r-\Xi\ \cot\theta\ \csc\theta}{{(r+\alpha \cos\theta)}^2}\partial_\phi+\frac{-r^3+\Xi(r+\alpha \cos\theta)+\alpha^2\ r}{{(r+\alpha \cos\theta)}^2}\partial_t \hfill
\end{multline}

\noindent Having obtained all Lie brackets for the elements of the vector frame, it is straightforward to show that Nijenhuis tensor vanishes for the almost complex structure by direct computation. We again use \eqref{Nijenhuis} for integrability check as in Taub-NUT case.

\noindent \begin{multline}
 N\left(e_2,e_4\right) = \left[e_2,e_4\right]+J \left[J \  e_2,e_4\right]+J \left[e_2,J \  e_4\right]-\left[J \  e_2 ,J \ e_4\right]  \quad \ \ \
\hfill\end{multline}~~\vspace{-6mm}
\begin{multline*}
\phantom{xxxxxx} = \left[e_2,e_4\right]+J \left[ e_3,e_4\right]-J \left[e_2, e_1\right]+\left[e_3 ,e_1\right]   \hfill
\end{multline*}~~\vspace{-6mm}
\begin{multline*}
\phantom{xxxxxx} = \frac{r\ \alpha\ \sin\theta \left((r^2-\alpha^2)\ \partial_t-\alpha\ \partial_\phi \right)}{\sqrt{\Delta} \ {(r+\alpha\ \cos\theta)}^2} + J\left(\frac{\alpha \sqrt{\Delta}\ (\cos\theta\ \partial_\theta -r \sin\theta\ \partial_r)}{{(r+\alpha\ \cos\theta)}^2} \right) \\ \phantom{xxxxxxxxxxx} -\frac{\alpha \sqrt{\Delta}\ \cot \theta\ (\partial_\phi + \alpha \sin^2 \theta\ \partial_t)}{{(r+\alpha\ \cos\theta)}^2}\hfill    
\end{multline*}~~\vspace{-6mm}
\begin{multline*}
\phantom{xxxxxx} =\frac{r\ \alpha\ \sin\theta \left((r^2-\alpha^2)\ \partial_t-\alpha\ \partial_\phi \right)}{\sqrt{\Delta} \ {(r+\alpha\ \cos\theta)}^2}  -\frac{\alpha \sqrt{\Delta}\ \cot \theta\ (\partial_\phi + \alpha \sin^2 \theta\ \partial_t)}{{(r+\alpha\ \cos\theta)}^2} \hfill
\end{multline*}
\begin{multline*}
\phantom{xxxxxx}+\frac{\alpha\sqrt\Delta \left (\cot\theta \ \partial_\phi-r\ \sin\theta \left(\partial_t+\frac{2M \ r \ \partial_t-\alpha \  \partial_\phi}{\Delta}\right)+\alpha \  \cos\theta \ \sin\theta \ \partial_t \right)}{{(r+\alpha \cos\theta)}^2} \\ \phantom{xxxxxxx} =0.
\hfill \end{multline*}

\noindent Again one example calculation is explicitly shown here, others can be easily computed by similar calculation steps. Note that $N\left(e_1,e_4\right)$ and $N\left(e_2,e_3\right)$ are almost trivially vanishing due to chosen almost complex structure. Thus, after conformally scaling the Euclidean Kerr metric in \eqref{Kerr2} we have shown that it becomes K{\"a}hler.
\,
\\

\noindent Now let's focus our attention to  alternative verification of the conformal K{\"a}hler structure. $W^+$ is needed in order to check the conditions of Theorem~\ref{Derdzinski}. It can be calculated for the  Euclidean Kerr metric in \eqref{Kerr2} by using the equation \eqref{Weyl+} as

\begin{equation} \label{Weyl2}
 A= W^+ + \frac{R}{12} I=\left(
\begin{array}{ccc}
 -\frac{M}{(r-\alpha  \cos \theta )^3} & 0 & 0 \\
 0 & -\frac{M}{(r-\alpha  \cos \theta )^3} & 0 \\
 0 & 0 & \frac{2 M}{(r-\alpha  \cos \theta)^3} \\
\end{array}
\right)   
\end{equation}

\noindent Kerr space is Ricci-flat ($R_{\mu\nu}=0$, $R=0$) so  $A=W^+ \neq 0$. It is 4-dimensional Einstein manifold and $W^+$ does not vanish, hence Theorem~\ref{Derdzinski} can be employed here. The matrix obtained in \eqref{Weyl2} indeed has only one distinct eigenvalue ($\lambda, \lambda, -2 \lambda$) and $|W^+|^2= \sum\limits_{i=1}^3 \lambda_i^2 $, where $\lambda_i$'s are the eigenvalues of  $W^+$. Then the conformal factor which transforms the original manifold to a K{\"a}hler manifold is computed as 
\begin{eqnarray}\label{Weylcon}
    |W^+|^{2/3}&=&\left(2\left( -\frac{M}{(r-\alpha  \cos \theta )^3}\right)^2+\left( \frac{2M}{(r-\alpha  \cos \theta )^3}\right)^2\right)^{1/3} \\ \nonumber
    &=&\frac{6^{1/3}M^{2/3}}{(r-\alpha  \cos \theta )^2}
\end{eqnarray}

\noindent Comparing the conformal factor in \eqref{Weylcon} with the one obtained from exact Lee-form, it is seen that they are same up to a constant. The arbitrary constant $C$ in \eqref{LeeKerr} becomes $6^{-1/3}M^{-2/3}$ for $ |W^+|^{2/3}$ factor . As mentioned before choice of arbitrary constant has no effect on the results. Therefore, we have verified the conformal K{\"a}hler structure from Weyl tensor analysis and found that the two conformal factors are matching.

\section{Conclusion}

In this paper, we presented alternative derivations of hyper-K{\"a}hler structure of Taub-NUT space and globally conformally K{\"a}hler structure of Kerr space via explicit coordinate-based calculations. We have obtained isometries between the most common explicit Taub-NUT metrics in literature and showed that Euclidean Kerr metric is globally conformally K{\"a}hler from two different analysis. We also showed that the conformal factor obtained from exact Lee-form matches with the one obtained from Weyl tensor analysis. This was shown for Kerr space, a more general proof may be needed but that can be subject of another study. It is our expectation that deriving K{\"a}hler structures of Kerr and Taub-NUT spaces starting directly from the commonly used metrics will provide insights for future studies which deal with complex geometry problems including these spaces.

\section*{Acknowledgments}

We thank M. Kalafat for useful discussions.

\end{document}